\documentclass[12pt]{article}
\usepackage{soul}
\usepackage[T1]{fontenc}
\usepackage[latin1]{inputenc}
\usepackage{bbm}
\usepackage{stmaryrd}
\usepackage{latexsym}
\usepackage{amsfonts}
\usepackage{amsmath,amssymb,amscd}
\usepackage{yfonts}
\usepackage{dsfont}
\usepackage{mathabx}
\usepackage[dvips]{graphicx}
\usepackage{epsfig}
\usepackage{psfrag}
\usepackage[font={small}]{caption}
\usepackage{color}
\usepackage{float}
\usepackage{upgreek}
\usepackage{tikz}
\usepackage{tikz-cd}
\usepackage{relsize}

\DeclareFontFamily{U}{txsyc}{}
\DeclareFontShape{U}{txsyc}{m}{n}{
   <-> txsyc%
}{}
\DeclareFontShape{U}{txsyc}{bx}{n}{
   <-> txbsyc%
}{}
\DeclareFontShape{U}{txsyc}{l}{n}{<->ssub * txsyc/m/n}{}
\DeclareFontShape{U}{txsyc}{b}{n}{<->ssub * txsyc/bx/n}{}
\DeclareSymbolFont{symbolsC}{U}{txsyc}{m}{n}
\SetSymbolFont{symbolsC}{bold}{U}{txsyc}{bx}{n}
\DeclareFontSubstitution{U}{txsyc}{m}{n}
\DeclareMathSymbol{\df}{\mathrel}{symbolsC}{"42}
\DeclareMathSymbol{\fd}{\mathrel}{symbolsC}{"43}
\DeclareMathSymbol{\lJoin}{\mathrel}{symbolsC}{"58}
\DeclareMathSymbol{\rJoin}{\mathrel}{symbolsC}{"59}

\newcommand{\cA}{{\cal A}}

\newcommand{\cE}{{\cal E}}

\newcommand{\cK}{{\cal K}}
\newcommand{\cL}{{\cal L}}

\newcommand{\cM}{{\cal M}}

\newcommand{\cP}{{\cal P}}

\newcommand{\cS}{{\cal S}}

\newcommand{\CC}{\mathbb{C}}

\newcommand{\LL}{\mathbb{L}}
\newcommand{\NN}{\mathbb{N}}
\newcommand{\PP}{\mathbb{P}}

\newcommand{\RR}{\mathbb{R}}

\newcommand{\ZZ}{\mathbb{Z}}

\newcommand{\fP}{\mathfrak{P}}

\newcommand{\fs}{\mathfrak{s}}

\newcommand{\iy}{\infty}

\newcommand{\lt}{\left}
\newcommand{\me}{\medskip}

\newcommand{\ri}{\rightarrow}
\newcommand{\rt}{\right}
\newcommand{\sm}{\smallskip}

\newcommand{\wi}{\widetilde}
\newcommand{\wit}{\widehat}

\newcommand{\card}{\mathrm{card}}

\DeclareMathOperator*{\esssup}{ess\,sup}

\newcommand{\fo}{\forall\ }

\newcommand{\lve}{\lt\vert}
\newcommand{\lVe}{\lt\Vert}

\newcommand{\rve}{\rt\vert}
\newcommand{\rVe}{\rt\Vert}

\newcommand{\stl }{\,:\,}

\newcommand{\un}{\mathds{1}}

\newcommand{\bq}{\begin{eqnarray*}}
\newcommand{\bqn}[1]{\begin{eqnarray}\label{#1}}
\newcommand{\eq}{\end{eqnarray*}}
\newcommand{\eqn}{\end{eqnarray}}
\newcommand{\wwtbp}{\par\hfill $\blacksquare$\par\me\noindent}
\newcommand{\thistitlepagestyle}{}
\newcommand{\lin}{\llbracket}
\newcommand{\rin}{\rrbracket}

\newcommand{\ttsim}{\raise.17ex\hbox{$\scriptstyle\mathtt{\sim}$}}

\newtheorem{pro}{Proposition}

\newtheorem{lem}[pro]{Lemma}

\renewcommand{\thepro}{\arabic{pro}}
\newenvironment{con}
{\par\me\refstepcounter{pro}\noindent{\bf Conjecture \thepro\ }}
{\par\hfill $\square$\par\sm\noindent}

\newenvironment{rem}
{\par\me\refstepcounter{pro}\noindent{\bf Remark \thepro\ }}
{\par\hfill $\square$\par\sm\noindent}

\newcommand{\proof}{\par\me\noindent\textbf{Proof}\par\sm\noindent}

\newcommand{\comment}[1]{}

\title{On intertwining relations between \\ Ehrenfest, Yule and Ornstein-Uhlenbeck processes}
\author{Laurent Miclo${}^\dagger$ and Pierre Patie${}^\ddagger$
}
 \date{\vbox{\copy0
 \vskip5mm
 \copy1
}
 }

\begin{document}

\setbox0=\vbox{
\large
\begin{center}
${}^\dagger$ Toulouse School of Economics, UMR 5314\\
Institut de Mathématiques de Toulouse, UMR 5219,\\
CNRS and Universit\'e de Toulouse
\end{center}
}

 \setbox1=\vbox{
 \large
 \begin{center}
 ${}^\ddagger$ School of Operations Research and Information Engineering,\\
 Cornell University \\
\end{center}
}
\setbox5=\vbox{
\hbox{${}^\dagger$ miclo@math.cnrs.fr\\}
\vskip1mm
\hbox{Toulouse School of Economics, UMR 5314,\\}
\hbox{Manufacture des Tabacs, 21, Allée de Brienne\\}
\hbox{31015 Toulouse cedex 6, France\\}
}
\setbox6=\vbox{
\hbox{${}^\ddagger$ pp396@cornell.edu\\}
\vskip1mm
\hbox{School of Operations Research and Information Engineering\\}
\hbox{Cornell University\\}
\hbox{Ithaca, NY 14853\\}
\hbox{USA\\}
 }

\maketitle
\thistitlepagestyle
\abstract{
Markovian intertwining relations between two Markov semigroups are related to the  partial inclusion of the spectra of their generators, at least for finite ergodic processes.
We check the limitations of this observation by investigating the Markov intertwining relations between the Ehrenfest, Yule and Ornstein-Uhlenbeck processes,
whose spectra are all included into $-\ZZ_+$. As a by-product, we offer a clarification of an intertwining relation found in Biane \cite{MR1459446} between the Yule and the Ornstein-Uhlenbeck processes.
}
\vfill\null
{\small
}\par

\newpage

\section{Introduction}

The state space reduction  is an important question in Markov process theory and its applications.
Given a Markov process $X\df(X_t)_{t\geq 0}$ on a large state space $V$, one is looking for another
Markov process $Y\df(Y_t)_{t\geq 0}$ on a much smaller state space $W$ and serving as a ``relatively good image'' of the evolution of $X$.
The process $Y$ corresponds to a limited quantity of information that one would like to extract from
$X$ while still providing for a sufficient knowledge about
 certain characteristics of the corresponding conditional distributions of the positions of $X$.
Some qualitative features are desirable in such an approximation/prediction procedure:
\begin{description}
\item[(i)]
The ``indicative process'' $Y$ takes into account the limited observation chosen to be made on $X$ in an non-anticipative way: for any $t\geq 0$,
to construct the piece of trajectory $Y_{[0,t]}$, we should only use what is extracted from $X$ up to time $t$ and maybe some additional independent randomness (which may require an enlargement of the underlying probability space, from a mathematical point of view). Furthermore, $Y_{[0,t]}$ is the only information we keep from our partial observations from $X_{[0,t]}$.
\item[(ii)]
The process $(X,Y)\df(X_t,Y_t)_{t\geq 0}$ is Markovian, to enable for ``online'' constructions. It is  time-homogeneous, as all the processes  considered here.
\item[(iii)] For any $t\geq 0$, knowing the trajectory $Y_{[0,t]}$, the conditional law of $X_t$ should depend only on $Y_t$,
to avoid the storage of too much information, since this is the objective of state space reduction. To be quite restrictive, we do not allow either for an explicit dependence on time.
\end{description}
Namely, we want to use  some partial observations of $X$ to construct in an adapted way a Markov process $Y$ whose current value $Y_t$ enables to make an ``as good as possible''
prediction on some aspects of the position $X_t$, given that we only observed $X$ through $Y$.
It may looks like filtering theory but it is different: there, the observation process $Y$ is given and we have to evaluate where is the signal process $X$.
Here we choose what to observe from $X$, encapsulated in $Y$, and it is limited because we want its state space to be small.\par\sm
Markov intertwinings meet the above requirements.
Initially they were developed by Diaconis and Fill \cite{MR1071805} in a discrete time and finite state space framework.
Let us recall the underlying principle in continuous time, as subsequently extended by Fill \cite{MR1144727}. The state spaces $V$ and $W$ are still assumed to finite and we are given $L^X$ the generator of $X$ on $V$.
In the first step, we look for a Markov generator $L^Y$ on $W$ and a Markov kernel $\Lambda$ from $W$ to $V$ such that
the following \textbf{intertwining relation} (said to go from $L^Y$ to $L^X$) holds
\bqn{interXY}
L^Y\Lambda&=&\Lambda L^X\eqn
Ideally, the Markov kernel $\Lambda$ should be the most ``informative'' possible, in particular its rank as a matrix should be $\min(\card(V),\card(W))$, which we expect to  be $\card(W)$ in the setting of state space reduction.
In the second step, when $Y$ is a Markov process generated by $L^Y$ and when its initial law $\cL(Y_0)$ satisfies $\cL(Y_0)\Lambda =\cL(X_0)$, we construct a coupling of
$X$ and $Y$ such that (i), (ii) and (iii) are satisfied:
\bqn{intertwining}
\fo t\geq 0,\qquad
\lt\{\begin{array}{rcl}
\cL(Y_{[0,t]}\vert X)&=&\cL(Y_{[0,t]}\vert X_{[0,t]})\\[2mm]
\cL((X,Y)_{[t,+\iy)}\vert (X,Y)_{[0,t]})&=&\cL((X,Y)_{[t,+\iy)}\vert (X_t,Y_t))\\[2mm]
 \cL(X_t\vert Y_{[0,t]})&=&\Lambda(Y_t, \cdot)
\end{array}\rt.\eqn
where the notation $\cL(\cdot\vert\cdot)$ stands for conditional laws.
\par\sm
To illustrate this procedure, let us come back to the historical example of the top-to-random shuffle due to Aldous and Diaconis \cite{MR841111}, in discrete time.
The state space is $V\df\cS_N$, the symmetric group on $N$ cards, and the transition of the Markov chain $X\df(X_n)_{n\in\ZZ_+}$
corresponds to taking the top card and replacing it at a uniformly chosen position in the deck of cards. Here we adopted the notation $\ZZ_+\df\{0, 1, 2, 3, ...\}$, while $\NN$ will stand for $\{1, 2, 3, ...\}$.
The Markov chain $Y\df(Y_n)_{n\in\ZZ_+}$ records the position of the card $C$ which initially was at the bottom of the deck, up to the time when it reaches the top of the deck
(when $C$ is replaced at random in the deck, by convention the current position of $Y$ is set at 0 and it stays there forever).  Thus $W=\lin 0,N\rin\df\{0,1, 2, ..., N\}$ and for large $N\in\NN$,
$\card(W)=N+1\ll N!=\card(V)$. In this example, we are interested in the distance in separation of the distribution of $X_n$ at time $n\in\ZZ_+$ to the invariant measure, which is the uniform distribution $\upsilon$ on $\cS_N$. Knowing the ``indicative process''  $(Y_m)_{m\in\lin 0, n\rin}$, we have a good idea of the distance in separation
of the conditional distribution of $X_n$ with $\upsilon$, in particular when $Y_n=0$, $\cL(X_n\vert (Y_m)_{m\in\lin 0, n\rin}, Y_n=0)=\upsilon$.
Nevertheless, this example does not convey very well the idea that for the purpose of state space reduction, we are rather looking for Markov kernels $\Lambda$ whose probability distributions $\Lambda(x,\cdot)$
do not spread much.
\par
\sm
Let us come back to the general situation. Given $L^X$ and $W,$ there is usually a lot of Markov generators $L^Y$ and Markov kernels $\Lambda$ such that \eqref{interXY} is satisfied.
So we must be a little more quantitative and wonder about what is a ``good intertwining relation''.
Note that if $\varphi$ is an eigenfunction associated to an eigenvalue $-\lambda\in\CC$ of $L^X$, then we get $L^Y[\psi]=-\lambda \psi$ with $\psi\df\Lambda[\varphi]$.
Namely, either $\psi=0$ or $\psi$ is an eigenfunction of $L^Y$ for the eigenvalue $-\lambda$.
\par
Conversely, assume that $V$ and $W$ are finite and that both the Markov generators
$L^X$ and $L^Y$ are irreducible. Suppose that some part $S$ of the spectrum of $L^Y$ is included into the spectrum of $L^X$.
Here spectrum has to be understood in an extended sense: it   concerns the size of the Jordan block as well as the value of the eigenvalue,
and multiplicity is taken into account. Then the computations of \cite{MR3838869}
enable to find
 a Markov kernel $\Lambda$ satisfying \eqref{interXY} and such that the image of $\Lambda$ contains the eigenspace for $L^Y$ associated to $S$.
There is a trivial instance of this principle: consider the case $S=\{0\}$, which is necessarily included into the spectra of $L^X$ and $L^Y$.
Then we can take for Markov kernel $\Lambda$ the invariant measure $\pi^X$ associated to $X$, namely we consider
\bq
\fo y\in W,\qquad \Lambda(y,\cdot)&=& \pi^X\eq
The general case is obtained by perturbation of this trivial situation.
In particular, $\Lambda$ may be quite small (measured for instance with respect to the image by $\Lambda$ of the unitary ball of $\LL^2(\pi^X)$) and a problem remains to find
the largest possible one.\par\sm
In the folklore,  when $\lambda$ is an eigenvalue of $L^X$, the smaller (respectively the larger)  is $\vert\lambda\vert$, the more $\lambda$ corresponds to global (resp.\ local)  features of the dynamics generated by $L^X$. For instance in the context of simulated annealing at small temperature, the smallest (non-zero) eigenvalue is directly related to
the largest height of a well not containing a fixed global minima of the underlying potential. If one wants to summarize such a process with a two-points dynamics, in some sense, one has to
cluster the well with the largest height into a unique point and its complementary set into the other point. This is the most global aspect of the full dynamics (after the fact that the process
does not lose mass, which corresponds to the eigenvalue zero). The following eigenvalues correspond to secondary features,
see for instance \cite{MR1364276} for their geometric description. Another instance of this heuristic is Weyl's law on a compact Riemannian manifold whose total volume is one (see e.g.\ the book of Taylor \cite{MR2743652}):
the behavior of the large eigenvalues of the Laplacian mainly depends on dimension of the manifold, since locally, manifolds of the same dimension all look  identical.
\par
The two above examples are somewhat asymptotical (one at small temperature and the other at large eigenvalues), nevertheless they suggest that if we are interested in the global behavior of the evolution of $X$,
we should rather look for intertwining relations \eqref{interXY} such that $\Lambda$ preserves the low lying part of the spectrum of $-L^X$ (while crushing the
eigenspaces corresponding to the remaining high lying part). The fact that the above examples are reversible (i.e.\ self-adjoint) is not relevant, it just insures that the  eigenvalues of the corresponding generators are real and non-positive. In general one has to consider the modules of the eigenvalues.
\par
\sm
These motivating observations lead us to the following problem. Given $L^X$ and a finite state space $W,$ find a Markov generator $L^Y$ and a Markov kernel $\Lambda$ from $W$ to $V$ such that \eqref{interXY} holds and the spectrum of $L^Y$ is the low lying spectrum of $L^X$.
In fact this is only the first part of the program described above, since  furthermore we would like $\Lambda$ to be the largest possible and also to couple the processes $X$ and $Y$.
The latter question is very important for applications, since given $X$, it amounts to knowing how to extract  the important information $Y$ from $X$.
In \cite{miclo:hal-01911989} we proposed a way to do it via the introduction of some random mappings in some particular situations where $Y$ is subset-valued (but then the state space of $Y$ can end up being much larger than the state space of $X$). \par
Here we will  only be concerned with a very special instance of this kind of issue, namely we will
consider some famous processes with the same low lying spectrum and we will try to find ``nice'' intertwining relations between them.
This is quite an academic point of view, but it will provide some preliminary insights on what it is possible to do and what is not, especially when the state space $V$ is infinite.\par\me
The first example we consider is the \textbf{Ehrenfest} family. For $N\in\ZZ_+$, define on $\lin 0, N\rin$ the  Markov generator $L_N$ via
\bqn{LN}
\fo x\not=x'\in \lin 0, N\rin,\qquad L_N(x,x')&\df&\frac12 \lt\{
\begin{array}{ll}
N-x&\hbox{, if $x'=x+1$}\\
x&\hbox{, if $x'=x-1$}\\
0&\hbox{, otherwise}
\end{array}\rt.
\eqn
(the values on the diagonal are such that the row sums all vanish).
It is well-known that the spectrum of $-L_N$ is $\lin 0, N\rin$.
\par\sm
The second example is the \textbf{Yule} family. For $N\in\ZZ_+$, consider on $\lin 0,N\rin$ the pure-death generator $D_N$ defined by
\bqn{DN}
\fo x\neq x'\in \lin 0, N\rin,\qquad D_N(x,x')&\df& \lt\{
\begin{array}{ll}
x&\hbox{, if $x'=x-1$}\\
0&\hbox{, otherwise}
\end{array}\rt.
\eqn\par
Since $D_N$ is a lower triangular matrix, its eigenvalues are given by the diagonal, namely $\lin 0, N\rin$ is the spectrum of $-D_N$.
Note that the generators of this family are not irreducible, as the associated processes are non-increasing.
There is a reverse family of Yule generators $(\wi D_N)_{N\in\ZZ_+}$ given by
\bqn{wiDN}
\fo N\in\ZZ_+,\,\fo x\neq x'\in \lin -N, 0\rin,\qquad \wi D_N(x,x')&\df& \lt\{
\begin{array}{ll}
-x+1&\hbox{, if $x'=x-1$}\\
0&\hbox{, otherwise}
\end{array}\rt.
\eqn\par
The spectrum of $-\wi D_N$ is also $\lin 0, N\rin$. We have already encountered these generators: up to a shift of the state space and a rescaling of time (after changing from discrete to continuous time), they correspond to the evolution of the last card in the top-to-random shuffle.\par
The family $( D_N)_{N\in\ZZ_+}$ admits an infinite version $D_{\iy}$: it is the pure-death  generator on $\ZZ_+$ whose infinite matrix $(D_{\iy}(y,y'))_{y,y'\in \ZZ_+}$ is imposed by its off-diagonal entries via:
\bq
\fo y\neq y'\in \ZZ_+,\qquad D_{\iy}(y,y')&\df& \lt\{
\begin{array}{ll}
y&\hbox{, if $y'=y-1$}\\
0&\hbox{, otherwise}
\end{array}\rt.
\eq\par
From a functional point of view, we see $D_{\iy}$ as an operator on $\CC^{\ZZ_+}$, via
\bq
\fo f\in \CC^{\ZZ_+},\,\fo y\in\ZZ_+,\qquad
D_{\iy}[f](y)&=& y(f(y-1)-f(y))\eq
\par
As it can be expected, the spectrum of $-D_{\iy}$ turns out to be $\ZZ_+$.\par
The reverse family $( \wi D_N)_{N\in\ZZ_+}$ equally admits an infinite version $\wi D_{\iy}$, on the state space $\ZZ_-\df-\ZZ_+$:
\bq
\fo y\neq y'\in \ZZ_-,\qquad \wi D_{\iy}(y,y')&\df& \lt\{
\begin{array}{ll}
1-y&\hbox{, if $y'=y-1$}\\
0&\hbox{, otherwise}
\end{array}\rt.
\eq
We will mainly work with the Yule family $(D_N)_{N\in\ZZ_+\sqcup\{\iy\}}$, since  Assertion (f) below reduces the interest of the reverse Yule family
(nevertheless, see the considerations at the beginning of Subsection  \ref{otrYf} and Conjecture \ref{conrev} at the end of this introduction).
\par
\sm
Our last Markov operator is the \textbf{Ornstein-Uhlenbeck} operator $L$ acting on $\cP$, the space of polynomial functions on $\RR$, via
 \bqn{OU}
 \fo f\in \cP,\,\fo x\in\RR,\qquad L[f](x)&\df&
 f''(x)-xf'(x)\eqn
 This operator is non-positive and symmetric in $\LL^2(\gamma)$, where $\gamma$ is the standard normal distribution.
Thus its Freidrichs extension provides  a self-adjoint operator on $\LL^2(\gamma)$, that is still  denoted by $L$.
 The spectrum of $-L$ consists of the eigenvalues  $n\in\ZZ_+$, all of them of multiplicity 1.
 \par\me
 The goal of this paper is to show the following assertions (under appropriate integrability assumptions for (b) and (d)):
 \begin{description}
 \item[(a)] There are surjective intertwinings from $L_N$ to $L_M$, for all $M\geq N\in\ZZ_+$.
  \item[(b)] The only intertwining from $L_N$ to $L$ is trivial, for all $ N\in\ZZ_+$.
  \item[(c)] There are surjective intertwinings from $D_N$ to $L_M$, for all $M\geq N\in\ZZ_+$.
  \item[(d)] The only intertwining from $L_N$ to $D_M$ is trivial, for any $ N,M\in\ZZ_+$.
   \item[(e)] There are surjective intertwinings from $D_N$ to $L$ for all $N\geq 2$, but not for $N=1$.
   \item[(f)] The only intertwining from $\wi D_N$ to $L$ is trivial for all $N\geq 1$.
   \end{description}
In these statements, an intertwining relation  is said to be trivial (respectively surjective) if the corresponding Markov kernel $\Lambda$ coincides with a probability distribution, i.e.\ if all its rows are the same (resp.\ if $\Lambda$ is surjective).
In (a) and (c), we will provide some explicit and quite natural  intertwinings. Concerning (e), we will describe all the possible intertwinings for $N=2$ and $N=3$, but for $N\geq 4$ the argument will only be perturbative, so some room is left for improvements that could lead to a proof of Conjecture \ref{con1} below.
The non-existence of (b) does not come from the fact we are trying to intertwine jump processes with diffusions, since in
\cite{2018arXiv180709445M} we intertwined infinite birth and death processes with Laguerre diffusions in non-trivial ways.\par
\sm
Although, as we shall after Lemma \ref{specdecomp}  below, the intertwining relation from $D_{\iy}$ to $L$ claimed in Theorem 3.4 of Biane \cite{MR1459446} does not hold, the original  abstract group theoretical approach developed in that paper  suggests that the following is true:
\begin{con}\label{con1}
There exists a non-trivial intertwining from $D_{\iy}$ to $L$.
\end{con}
Despite (f), we equally believe in:
\begin{con}\label{conrev}
There exists a non-trivial intertwining from $\wi D_{\iy}$ to $L$.
\end{con}
\par
In the next section we study  the intertwining relations starting from an ergodic generator, namely (a), (b) and (d).
Section \ref{isfaag} deals with the remaining cases, where the intertwining relations starts from an absorbed generator.

 \bigskip
 \par\hskip5mm\textbf{\large Acknowledgments:}\par\sm\noindent
The authors  are grateful to the Toulouse School of Economics  where this work was done and initiated during the invitation of the second author.

\section{Intertwinings
from an ergodic generator}

Here we deal with the points (a), (b) and (d) of the introduction, respectively in the following subsections.

\subsection{From Ehrenfest to Ehrenfest}\label{fEtE}

Fix some $N\in\ZZ_+$. The invariant probability measure $\pi_N$ associated to the generator $L_N$ defined in \eqref{LN} is the binomial distribution given by
\bqn{piN}
\fo x\in \lin 0, N\rin,\qquad \pi_N(x)&\df& 2^{-N}\binom{N}{x}\eqn
\par
This measure  is furthermore reversible for $L_N$, i.e.\
the operator $L_N$ is self-adjoint in $\LL^2(\pi_N)$.
It follows that $L_N$ is diagonalizable. It is well-known that the set of eigenvalues of $-L_N$ is $\lin 0, N\rin$, all with multiplicity 1.
The eigenspace associated to the eigenvalue $n\in \lin 0, N\rin$ is generated by the Krawtchouk polynomial $K_{N,n}$.
These polynomials can be defined via their generating function:
\bqn{GN}
 \nonumber
 \fo x\in \lin 0, N\rin,\,\fo z\in \CC,\qquad
 G_N(z,x)&\df&
\sum_{n\in \lin 0, N\rin} K_{N,n}(x)\frac{z^n}{n!}\\
&=&\lt( 1+\frac{z}{2}\rt)^x\lt(1-\frac{z}{2}\rt)^{N-x}\eqn
see for instance Griffiths \cite{MR3501191}.
\par
Since the set of eigenvalues of $L_N$ is included into the set of eigenvalues of $L_{N+1}$ and the probability measure $\pi_{N+1}$ charges all the points of $\lin 0, N+1\rin$, from \cite{MR3838869}, we get  that there exists a Markov kernel $\Lambda_N$ from $\lin 0, N\rin$ to $\lin 0, N+1\rin$
such that
\bqn{interN}
L_N\Lambda_N&=&\Lambda_N L_{N+1}\eqn
and such that the rank of $\Lambda_N$, seen as a $\lin 0, N\rin\times \lin 0, N+1\rin$-matrix, is $N+1$.
\par
Let us exhibit a very simple one:
\begin{lem}\label{NN}
Consider the Markov kernel $\Lambda_N$ from $\lin 0, N\rin$ to $\lin 0, N+1\rin$ given by
\bq
\fo x\in \lin 0, N\rin,\,\fo y\in \lin 0, N+1\rin,\qquad \Lambda_N(x,y)&\df&
 \lt\{
\begin{array}{ll}
1/2&\hbox{, if $y=x+1$ or $y=x$}\\
0&\hbox{, otherwise}
\end{array}\rt.
\eq
We have
\bqn{LKKN2}
\fo n\in \lin 0, N+1\rin,\qquad \Lambda_N[K_{N+1,n}]&=& K_{N,n}\eqn
with the convention that $K_{N,N+1}= 0$ on $\lin 0, N\rin$. In particular the intertwining relation \eqref{interN} is satisfied.
\end{lem}
\proof
We compute the generating function of the family $(\Lambda_N[K_{N+1,n}])_{n\in \lin 0, N+1\rin}$: for any $z\in \CC$ and $x\in \lin 0, N\rin$,
\bq
 \sum_{n\in \lin 0, N+1\rin} \Lambda[K_{N+1,n}](x)\frac{z^n}{n!}
 &=&\frac12\sum_{n\in \lin 0, N+1\rin} (K_{N+1,n}(x)+K_{N+1,n}(x+1))\frac{z^n}{n!}\\
&=&\frac12\lt(\sum_{n\in \lin 0, N+1\rin} K_{N+1,n}(x)\frac{z^n}{n!}+\sum_{n\in \lin 0, N+1\rin} K_{N+1,n}(x+1)\frac{z^n}{n!}\rt)\\
&=&\frac12\lt( \lt( 1+\frac{z}{2}\rt)^x\lt(1-\frac{z}{2}\rt)^{N+1-x}+\lt( 1+\frac{z}{2}\rt)^{x+1}\lt(1-\frac{z}{2}\rt)^{N-x}\rt)\\
&=&\frac12 \lt(1+\frac{z}{2}\rt)^x\lt(1-\frac{z}{2}\rt)^{N-x}\lt(1-\frac{z}{2}+ 1+\frac{z}{2}\rt)\\
&=& \lt(1+\frac{z}{2}\rt)^x\lt(1-\frac{z}{2}\rt)^{N-x}
\eq
where \eqref{GN} was taken into account in the second equality.
Using again \eqref{GN}, we deduce
\eqref{LKKN2}.
To get the sought intertwining relation, it is sufficient to check it on the $(K_{N+1,n})_{n\in \lin 0, N+1\rin}$,
 which is a base of $\LL^2(\pi_{N+1})$. We have for any $n\in \lin 0, N+1\rin$,
 \bq
 L_N\Lambda_N[K_{N+1,n}]&=&L_N [K_{N,n}]\\
 &=&-n K_{N,n}\\
 &=&-n \Lambda_N[K_{N+1,n}]\\
 &=&\Lambda_NL_{N+1}[K_{N+1,n}]\eq
\wwtbp
\par
The relation \eqref{interN} can be extended to other pairs of Ehrenfest generators. Indeed,
writing 
\bq
\fo M\leq N\in \ZZ_+,\qquad
\Lambda_{M,N}&\df& \Lambda_{M}\Lambda_{M+1}\cdots \Lambda_{N-1}\eq
we get  that for any $ M\leq N\in \ZZ_+$,
\bqn{LNMxy}
\fo x\in \lin 0, M\rin,\fo y\in \lin 0, N\rin,\qquad \Lambda_{M,N}(x,y)&=&
\lt\{
\begin{array}{ll}
2^{M-N} \binom{N-M}{y-x}&\hbox{, if $x\leq y\leq x+N-M$}\\
0&\hbox{, otherwise}
\end{array}\rt.
\eqn
\par
By an immediate iteration of \eqref{LKKN2}, we obtain the intertwining relation
\bqn{LNLNM2}
\fo M\leq N\in \ZZ_+,\qquad L_M\Lambda_{M,N}&=&\Lambda_{M,N} L_{N}\eqn

\subsection{From Ehrenfest to Ornstein-Uhlenbeck}

 Consider the family of Hermite polynomials $(h_n)_{n\in\ZZ_+}$, defined, similarly to \eqref{GN}, via their generating function:
 \bqn{hermite}
  \fo x\in \RR,\,\fo z\in \CC,\qquad
\sum_{n\in \ZZ_+} h_n(x)\frac{z^n}{n!}
&=&\exp(zx-z^2/2)\eqn
For any $n\in\ZZ_+$, the eigenspace associated to the eigenvalue $-n$ of the Ornstein-Uhlenbeck generator $L$ defined in \eqref{OU} is
generated by $h_n$.
\par
Fix some $N\in\ZZ_+$. We are wondering whether there exists a Markov kernel $\Lambda$ from $\lin 0, N\rin$ to $\RR$ such that
\bqn{LNLLL}
L_N\Lambda&=&\Lambda L\eqn
This equality is understood as holding on $\cP$, so it is implicitly assumed that for any $n\in\ZZ_+$, the probability measure $\Lambda(n,\cdot)$ has moments of all orders.\par
Since the Gaussian measure $\gamma$ is invariant for $L$, \eqref{LNLLL} is satisfied with $\Lambda=\gamma$, namely with the trivial Markov kernel
given by
\bqn{trivial}
\fo n\in\ZZ_+,\qquad \Lambda(n,\cdot)&=& \gamma\eqn
Indeed, with this Markov kernel, both sides of \eqref{LNLLL} vanish.
\par
Our goal here is to show that under a strengthened integrability assumption, $\gamma$ is the unique Markov kernel such that \eqref{LNLLL} is satisfied.
\par
Let $\cM$ be the set of probability measures on $\RR$ integrating the mapping $\RR\ni x\mapsto \exp(x^2/4)$ and
denote by $\cK_N$ the set of Markov kernels $\Lambda$  from $\lin 0,N\rin$ to $\RR$
 such that $\Lambda(n,\cdot)\in\cM$ for all $n\in \lin 0, N\rin$.
We have:
\begin{pro}\label{EhOU}
The only Markov kernel $\Lambda\in\cK_N$ such that \eqref{LNLLL} is satisfied is the trivial kernel defined in \eqref{trivial}.
\end{pro}
\proof
For any $n\in\ZZ_+$, consider $\psi_n\df \Lambda[h_n]$.
Due to \eqref{LNLLL}, we have
\bq
L_N[\psi_n]&=&L_N\Lambda[h_n]\\
&=&\Lambda L[h_n]\\
&=&-n\Lambda[h_n]\\
&=&-n\psi_n\eq
It follows that $\psi_n$ is an eigenvector associated to the eigenvalue $-n$ of $L_N$.
According to Subsection~\ref{fEtE}, $\psi_n$ is proportional to the Krawtchouk polynomial $K_{N,n}$, with the convention that $K_{N,n}=0$ on $\lin 0, N\rin$ for $n>N$.
Let $a_n\in\RR$ be such that $\psi_n=a_nK_{N,n}$, in particular $a_n=0$ for $n>N$.\par
With the help of \eqref{hermite} and Lemma \ref{Fubini} below, which justifies the application of Fubini's lemma, we compute that for any $z\in\CC$ and $y\in \lin 0, N\rin$,
\bqn{erreur3}
\int_{\RR} \Lambda(y,dx) \exp(zx-z^2/2)&=&
\sum_{n\in\ZZ_+} \Lambda[h_n](y) \frac{z^n}{n!}\\
\nonumber&=&\sum_{n\in\ZZ_+} a_n \frac{z^n}{n!} K_{N,n}(y)\\
\label{erreur4}&=&\sum_{n\in\lin 0, N\rin} a_n \frac{z^n}{n!} K_{N,n}(y)
\eqn
\par
Conversely, observe that if $\Lambda$ is a Markov kernel from $\ZZ_+$ to $\RR$ satisfying \eqref{erreur4} for all $z\in\CC$ and $y\in \lin 0, N\rin$,
then \eqref{LNLLL} is true. Indeed, from \eqref{erreur4} we deduce that $\psi_n=a_nK_{N,n}$ for all $n\in\ZZ_+$,
namely $\Lambda$ maps each eigenspace of $L$ into the corresponding eigenspace of $L_N$ (with the convention that the eigenspaces associated to the eigenvalues $-n$ with $n>N$ are reduced to $\{0\}$) and this is sufficient to insure \eqref{LNLLL}.\par\sm
 Taking into account \eqref{hermite} and the fact  that $(h_n)_{n\in\ZZ_+}$ is an orthogonal family in $\LL^2(\gamma)$, we get that for any $n\in\ZZ_+$,
 \bqn{hgz}
 \int h_n(x) \exp(zx-z^2/2)\, \gamma(dx) &=& z^n\eqn
 due to the fact that
 \bq
\fo n\in\ZZ_+, \qquad \int h_n^2(x) \, \gamma(dx) &=&n!\eq\par
It follows that for any $y\in\lin 0, N\rin$,
\bq
\int_{\RR} \lt(\sum_{n\in\lin 0, N\rin}  K_{N,n}(y)\frac{a_n}{n!} h_n(x)\rt) \exp(zx)\,\gamma(dx) &=&\lt(\sum_{n\in\lin 0, N\rin} K_{N,n}(y)\frac{a_n}{n!}z^n \rt)\exp(z^2/2)\\
&=&\int_{\RR} \Lambda(y,dx) \exp(zx)
\eq
We deduce that for all $y\in\lin 0, N\rin$, we have
\bqn{pasM}
\Lambda(y,dx)&=&\sum_{n\in\lin 0, N\rin}  K_{N,n}(y)\frac{a_n}{n!} h_n(x)\gamma(dx)\eqn
In particular, we must have for all $y\in\lin 0, N\rin$
and a.e.\ $x\in\RR$,
\bq
\sum_{n\in\lin 0, N\rin}  K_{N,n}(y)\frac{a_n}{n!} h_n(x)&\geq &0\eq
By continuity in $x$ of the left-hand side, this should be true for all $y\in\lin 0, N\rin$ and $x\in\RR$.
\par
Let $n_0\in \lin 0, N\rin$ be the largest integer such that $a_{n_0}\neq 0$. Assume that $n_0\geq 1$.
Since for any $n\in\ZZ_+$, the polynomial $h_n$ has degree $n$ and its highest coefficient is 1,
we have that as $x$ goes to $+\iy$,
\bq
\sum_{n\in\lin 0, N\rin}  K_{N,n}(y)\frac{a_n}{n!} h_n(x)&\sim& K_{N,n_0}(y)\frac{a_{n_0}}{n_0!} x^{n_0}\eq
Since the integral of  $K_{N,n_0}$  with respect to $\pi_N$ is zero (the scalar product in $\LL^2(\pi_N)$ of $K_{N,n_0}$ and $K_{N,0}=\un$, the function always taking the value 1, vanishes),
we can find $y_-, y_+\in \lin 0, N\rin$ such that $K_{N,n_0}(y_-)<0$ and $K_{N,n_0}(y_+)>0$.
Thus we can find $y_0\in\{y_-, y_+\}$ such that $K_{N,n_0}(y_0)\frac{a_{n_0}}{n_0!}<0$, namely
\bq
\lim_{x\ri+\iy} \sum_{n\in\lin 0, N\rin}  K_{N,n}(y_0)\frac{a_n}{n!} h_n(x)&=&-\iy\eq
This is in contradiction with the non-negativity of the left-hand-side, so necessarily $n_0=0$.
We deduce that
for all $y\in\lin 0, N\rin$, we have
\bq
\Lambda(y,dx)&=& K_{N,0}(y)\frac{a_0}{0!} h_0(x)\gamma(dx)\\
&=&a_0 \gamma(dx)
\eq
We get that $a_0=1$ and finally $\Lambda=\gamma$.
\wwtbp
The above proof is not yet complete, since we did not validate the use of Fubini's lemma in \eqref{erreur3}. This is done in the next result, which will also justify the integrability assumption of Proposition \ref{EhOU}.
\begin{lem}\label{Fubini}
For any given $y\in \ZZ_+$, the identity \eqref{erreur3} is justified as soon as $\Lambda(y,\cdot)\in\cM$.
\end{lem}
\proof
It is a consequence of Cramer's inequality, see for instance the book of Szeg\"{o} \cite{MR0372517}, claiming that there exists a constant $c>0$ such that
\bq
\fo n\in\ZZ_+,\,\fo x\in\RR,\qquad \vert h_n(x)\vert &\leq & c\sqrt{n!}\exp(x^2/2)\eq
Indeed, this bound yields  that for any given $z\in\CC$,
\bq
\int_{\RR} \Lambda(y,dx) \sum_{n\in\ZZ_+}\lve h_n(x)\frac{z^n}{n!}\rve& \leq &
c\int_{\RR} \Lambda(y,dx) \sum_{n\in\ZZ_+}\frac{\exp(x^2/4)\vert z\vert^n}{\sqrt{n!}}\\
&=&c\int_{\RR} \Lambda(y,dx) \exp(x^2/4)\sum_{n\in\ZZ_+}\frac{\vert z\vert^n}{\sqrt{n!}}\\
&<&+\iy\eq
\wwtbp\par
\begin{rem}\label{pasM2}
If one is not interested in Markov kernels, note that the signed kernels defined in \eqref{pasM} always provide  intertwining links between Ehrenfest and Ornstein-Uhlenbeck generators.
\end{rem}
Proposition \ref{EhOU} is an unpleasant fact for the program described in the introduction: consider any irreducible Markov generator $G_N$ on $\lin 0, N\rin$ whose
eigenvalues are $-\lin 0, N\rin$, where $N\in\ZZ_+$ is fixed. Then there is no non-trivial intertwining
from $G_N$ to $L$, namely it is not possible to have a good image (in the sense of the introduction) on $N+1$ points of the Ornstein-Uhlenbeck process.
Indeed, assume on the contrary that we have such an non-trivial intertwining $G_N\Lambda=\Lambda L$.
Since $L_N$ and $G_N$ have the same spectrum and are irreducible, we deduce from \cite{MR3838869}  there exists a Markov kernel $\wi\Lambda$ from $\lin 0,N\rin$ into itself such that
$L_N\wi\Lambda=\wi\Lambda G_N$ and such that $\wi\Lambda$ is an invertible matrix.
We would then obtain the non-trivial intertwining relation $L_N\wi\Lambda\Lambda=\wi\Lambda\Lambda L$, a contradiction.

\subsection{From Ehrenfest to Yule}

In general the Markov kernel entering into an intertwining relation transports the invariant measure of the first generator into the invariant measure of the second generator.
So when the second generator is absorbing with a unique absorbing point, the Markov kernel must be trivial and equal to the Dirac mass on the absorbing point.\par
\sm
Let us illustrate this principle on intertwining relations between Ehrenfest and Yule processes. For any given $N,M\in\ZZ_+$, assume that $\Lambda$ is a Markov kernel from $\lin 0, N\rin$ to $\lin 0, M\rin$ such that
\bqn{LNLLDM}
L_N\Lambda&=&\Lambda D_M\eqn
Let $\pi_N$, described in \eqref{piN}, be the invariant probability of $L_N$, for instance seen as a row vector.
Multiplying \eqref{LNLLDM} on the left by $\pi_N$, we get
\bq
\pi_N\Lambda D_M&=&\pi_NL_N\Lambda\ =\ 0\eq
Thus the probability $\pi_N\Lambda$ on $\lin 0, M\rin$ is an invariant probability for $D_M$. Since the Markov processes associated to $D_M$ all end up being absorbed at 0,
necessarily $\pi_N\Lambda=\delta_0$. Taking into account that $\pi_N$ charges all the points of $\lin 0, N\rin$, we get that
\bq
\fo x\in \lin 0, N\rin,\qquad \Lambda(x,\cdot)&=&\delta_0\eq
namely $\Lambda$ is trivial.\par
\sm
The above arguments are immediately extended to the case where $M=\iy$ and also to intertwining relations of the form $L\Lambda=\Lambda D_M$, with $M\in\ZZ_+\sqcup\{\iy\}$.
We conclude that a.e.\ in $x\in\RR$, $\Lambda(x,\cdot)=\delta_0$.
However for $L\Lambda=\Lambda D_M$ to make sense a priori, we must assume that $\Lambda$ transforms $\RR^{\lin 0, M\rin}$ into functions that are at least continuous (or seen as elements of $\LL^2(\gamma)$, if we consider the Friedrich extension), so that $\Lambda(x,\cdot)=\delta_0$ holds for all $x\in\RR$ (or $\Lambda=\delta_0$ in the $\LL^2(\gamma)$ context).

\section{Intertwinings
from an absorbed generator}\label{isfaag}

In the last subsection we have seen there is usually no non-trivial intertwining relation from an ergodic generator to an absorbed generator.
It is not true in the reverse direction, as shown by the top-to-random card shuffle example of Aldous and Diaconis \cite{MR841111}, and more generally such intertwining relations were exploited by Diaconis and Fill \cite{MR1071805}
to construct strong stationary times.
Here we deal with the points (c), (e) and (f) of the introduction, respectively in the following subsections.

\subsection{From Yule to Ehrenfest}

The intertwining relation from Yule to Ehrenfest is well-known, as well as the relations with the discrete hypercube random walks, see Example 4.38 of Diaconis and Fill \cite{MR1071805}.
We give here a direct proof in the spirit of this paper, via generating functions.\par\sm
Fix some $N\in\ZZ_+$. We have seen that the eigenvalues of the Yule generator $D_N$ defined in \eqref{DN} are the elements of $-\lin 0, N\rin$.
Let us compute the corresponding eigenvectors:
\begin{lem}\label{varphin}
For all $n\in \lin 0, N\rin$, the eigenspace associated to the eigenvalue $n$ of $-D_N$ is generated by the function $\varphi_n$ given by
\bqn{defphi}
\fo x\in \lin 0, N\rin,\qquad \varphi_n(x)&\df&
\binom{x}{n}\eqn
with the usual convention that $\binom{m}{n}=0$ for any $m< n$.
\end{lem}
\proof
Fix $n\in \lin 0, N\rin$ and let us check that for any $x\in \lin 0, N\rin$,
\bqn{DNvarphi} D_N[\varphi_n](x)&=&-n\varphi_n(x)\eqn
\par
For $x\in \lin 0, n-1\rin$, both sides are zero, so the equality holds.
\par
For $x=n$, the left-hand side is equal to $n(\varphi_n(n-1)-\varphi_n(n))=-n\varphi_n(n)$,
so the equality holds.\par
Assume that \eqref{DNvarphi} is true for some $x\geq n$ and let us show it is also true with $x$ replaced by $x+1$, as long as $x+1\leq N$.
We have
\bq
D_N[\varphi_n](x+1)&=&(x+1)(\varphi_n(x)-\varphi_n(x+1))\\
&=&(x+1)\lt(\binom{x}{n}-\binom{x+1}{n}\rt)\\
&=&\frac{x+1}{n!}\big( x(x-1)\cdots (x-n+1)-(x+1)x\cdots (x-n+2)\big)\\
&=&\frac{x+1}{n!}x(x-1)\cdots (x-n+2)\big(x-n+1-(x+1)\big)\\
&=&-n\frac{x+1}{n!}x(x-1)\cdots (x-n+2)\\
&=&-n\varphi_n(x+1)
\eq
which completes the proof.
\wwtbp
\par
We proceed by computing  the  generating function of the family $(2^{-n}n!\varphi_n)_{n\in \lin 0, N\rin}$:
\begin{lem}
We have
\bq
\fo z\in\CC,\,\fo x\in \lin 0, N\rin,\qquad
\sum_{n\in\lin 0,N\rin} 2^{-n}n!\varphi_n(x)\frac{z^n}{n!}&=& \lt(1+\frac{z}{2}\rt)^x\eq
\end{lem}
\proof
Indeed, we compute that
\bq
\sum_{n\in\lin 0,N\rin} 2^{-n}n!\varphi_n(x)\frac{z^n}{n!}&=&\sum_{n\in\lin 0,N\rin} \varphi_n(x)(2^{-1}z)^n\\
&=&\sum_{n\in\lin 0,N\rin}\binom{x}{n}(2^{-1}z)^n\\
&=&\lt(1+\frac{z}{2}\rt)^x\eq
\wwtbp
\par
Consider the Markov kernel $\wit \Lambda_N$ from $\lin 0, N\rin$ to $\lin 0, N\rin$ defined by
\bq
\fo x,y\in \lin 0, N\rin,\qquad
\wit \Lambda_N(x,y)&\df&2^{x-N}\binom{N-x}{y-x}\\
&=& \Lambda_{x,N}(x,y)
\eq
with the notation introduced in \eqref{LNMxy}.
In particular the support of $\Lambda(x,\cdot)$ is $\lin x, N\rin$.\par
Its introduction is motivated by:
\begin{pro}
We have
\bq
\fo n\in \lin 0, N\rin,\qquad \wit \Lambda_N[K_{N,n}]&=&2^{-n}n!\varphi_n\eq
and in particular,
\bqn{DNLN}
D_N\wit \Lambda_N&=&\wit \Lambda_N L_N\eqn
\end{pro}
\proof
For the first equality, it is sufficient to check the equality of the generating functions associated to the families $(\wit \Lambda_N[K_{N,n}])_{n\in \lin 0, N\rin}$
and
 $(2^{-n}n!\varphi_n)_{n\in \lin 0, N\rin}$, namely
 \bq
 \fo z\in\CC,\,\fo x\in \lin 0, N\rin,\qquad \sum_{n\in \lin 0, N\rin} \wit \Lambda_N[K_{N,n}](x)\frac{z^n}{n!}
 &=& \lt(1+\frac{z}{2}\rt)^x\eq
With the notation of \eqref{GN}, this equality is equivalent to
 \bqn{con2}
 \fo z\in\CC,\,\fo x\in \lin 0, N\rin,\qquad  \wit \Lambda_N[G_{N}(z,\cdot)](x)
 &=& \lt(1+\frac{z}{2}\rt)^x\eqn
so let us compute the l.h.s.: for any  $z\in\CC$ and $x\in \lin 0, N\rin$,
 \bq
 \wit \Lambda_N[G_{N}(z,\cdot)](x)&=&2^{x-N}\sum_{y\in \lin x, N\rin}\binom{N-x}{y-x}G_N(z,y)\\
 &=&2^{x-N}\sum_{y\in \lin x, N\rin}\binom{N-x}{y-x} \lt( 1+\frac{z}{2}\rt)^y\lt(1-\frac{z}{2}\rt)^{N-y}\\
 &=&2^{x-N} \lt( 1+\frac{z}{2}\rt)^x\sum_{y\in \lin x, N\rin}\binom{N-x}{y-x} \lt( 1+\frac{z}{2}\rt)^{y-x}\lt(1-\frac{z}{2}\rt)^{N-y}\\
 &=&2^{x-N} \lt( 1+\frac{z}{2}\rt)^x\lt(1+\frac{z}{2}+1-\frac{z}{2}\rt)^{N-x}\\
&=& \lt( 1+\frac{z}{2}\rt)^x
 \eq
This enables us to conclude the validity of \eqref{con2}.\par
Concerning the intertwining equality announced in the lemma, it is sufficient to check it on the basis $(K_{N,n})_{n\in \lin 0, N\rin}$ of $\RR^{\lin 0, N\rin}$.
Indeed, for any $n\in \lin 0, N\rin$, we have
\bq
\wit \Lambda_N L_N[K_{N,n}]&=&\Lambda[-n K_{N,n}]\\
&=&-n\Lambda[K_{N,n}]\\
&=&-n 2^{-n}n!\varphi_n\\
&=&D_N[2^{-n}n!\varphi_n]\\
&=&D_N\wit \Lambda_N [K_{N,n}]\eq
\wwtbp
\par\sm
Note that the kernel $\wit \Lambda_N$ is surjective. For integers $M\leq N$, one deduces a surjective intertwining relation
\bqn{DNLMN}
D_M\wit\Lambda_{M,N}&=&\wit \Lambda_{M,N} L_N\eqn
from \eqref{DNLN}
 by restriction from $\lin 0, N\rin$ to $\lin 0, M\rin$.
Formally, consider $I_{M,N}$ is the natural imbedding from  $\lin 0, M\rin$ into $\lin 0, N\rin$, seen as a Markov kernel.
Since $D_N$ is a lower triangular matrix, we have the simple intertwining relation $I_{M,N}D_N=D_MI_{M,N}$.
Thus multiplying \eqref{DNLN} on the left by $I_{M,N}$, we obtain \eqref{DNLMN} with the surjective kernel
$\wit \Lambda_{M,N}\df I_{M,N}\wit \Lambda_N$.
\par
Alternatively, we can start from $D_M\wit\Lambda_M=\wit\Lambda_M L_M$ that we multiply on the right by the Markov kernel $\Lambda_{M,N}$ defined in \eqref{LNMxy}.
Taking into account \eqref{LNLNM2}, we end up with \eqref{DNLMN} with
$\wit \Lambda_{M,N} \df \wit\Lambda_N\Lambda_{M,N}$. It is not difficult to see that the latter Markov kernel coincides with the former one, i.e.\ $\wit\Lambda_M\Lambda_{M,N}=I_{M,N}\wit \Lambda_N$. Indeed, it is sufficient to check this identity when it is applied to the family of Krawtchouk polynomials $(K_{N,n})_{n\in\lin 0, N\rin}$.

\subsection{From Yule to Ornstein-Uhlenbeck}\label{fYtO}

We start by considering the finite Yule generators and thus fix some $N\in\ZZ_+$.
Recall that the set $\cK_N$ of Markov kernels was defined before Proposition \ref{EhOU}. We are interested in
\bq
\cL_N&\df&\{\Lambda \in\cK_N\stl  D_N\Lambda\ =\ \Lambda L\}\eq
\par
To any $a\df(a_n)_{n\in\lin 0,N\rin}$, we associate the mapping $\lambda_a\stl  \lin0,N\rin\times\RR\ri\RR$
and the signed kernel $\Lambda_a$ from $\lin 0, N\rin$ to $\RR$ via
\bqn{lama}
\fo y\in\lin 0, N\rin,\,\fo x\in \RR,\qquad \lambda_a(y,x)&\df& \sum_{n\in\lin 0,y\rin} \frac{a_n}{n!}\binom{y}{n}h_n(x)\\
\label{Lama} \fo y\in\lin 0, N\rin,\qquad \Lambda_a(y,dx)&\df& \lambda_a(y,x)\, \gamma(dx)
\eqn
Finally, consider $\cA_N$ the set of $a\df(a_n)_{n\in\lin 0,N\rin}\in \RR^{N+1}$ with $a_0=1$ such that
\bqn{cAN}
\fo y\in\lin 0, N\rin,\,\fo x\in \RR,\qquad \lambda_a(y,x)&\geq &0\eqn
 \par
 The interest of these definitions is:
 \begin{pro}\label{comput}
 We have
 \bq
 \cL_N&=&\{\Lambda_a\stl  a\in\cA_N\}\eq
 \end{pro}
 As in Remark \ref{pasM2}, the following proof shows that the signed kernels $\Lambda_a$ always provide intertwining links between the finite Yule and Ornstein-Uhlenbeck generators.
 \proof
 The arguments are similar to those of the proof of Proposition \ref{EhOU}.
 Consider $\Lambda\in\cL_N$ and for $n\in\lin 0, N\rin$ define $\psi_n\df\Lambda[h_n]$.
 The intertwining relation $D_N\Lambda = \Lambda L$ implies that $D_N[\psi_n]=-n\psi_n$
 and according to Lemma \ref{varphin}, $\psi_n$ is proportional to $\varphi_n$.
 Denote $a_n\in\RR$ such that $\psi_n=a_n\varphi_n$.
 Note that when $n=0$, we have $\psi_0=\un=\varphi_0$, so that $a_0=1$.
 \par
 With the help of \eqref{hermite} and Lemma \ref{Fubini}, we get that for any $z\in\CC$ and $y\in \ZZ_+$,
\bqn{erreur2}
\nonumber\int_{\RR} \Lambda(y,dx) \exp(zx-z^2/2)&=&
\sum_{n\in\ZZ_+} \Lambda[h_n](y) \frac{z^n}{n!}\\
\nonumber&=&\sum_{n\in\ZZ_+} a_n \frac{z^n}{n!} \varphi_n(y)\\
\nonumber&=&\sum_{n\in\ZZ_+} a_n \frac{z^n}{n!} \binom{y}{n}\\
&=&\sum_{n\in\lin 0,y\rin} a_n \frac{z^n}{n!} \binom{y}{n}
\eqn
\par
Conversely, observe that if $\Lambda$ is a Markov kernel from $\ZZ_+$ to $\RR$ satisfying \eqref{erreur2} for all $z\in\CC$ and $y\in \ZZ_+$,
then the intertwining relation $D_N\Lambda = \Lambda L$ holds. Indeed, from \eqref{erreur2} we deduce that $\psi_n=a_n\varphi_n$ for all $n\in\lin 0,N\rin$,
namely $\Lambda$ maps each eigenspace of $L$ into the corresponding eigenspace of $D$
(again with the convention that the eigenspaces associated to the eigenvalues $-n$ with $n>N$ are reduced to $\{0\}$)
and this is sufficient to insure that $D_N\Lambda=\Lambda L$.\par\sm
Taking into account \eqref{hgz}, it appears that for any $z\in\CC$ and $y\in \ZZ_+$,
\bq
\int_{\RR} \Lambda(y,dx) \exp(zx)&=&\int_{\RR} \Lambda_a(y,dx) \exp(zx)\eq
with $a\df(a_n)_{n\in\lin 0, N\rin}$, namely $\Lambda=\Lambda_a$.
Since $\Lambda$ is a non-negative kernel, we get that $a\in\cA_N$. Conversely,
when $a\in \cA_N$, it remains to check that
the total mass of $\Lambda_a(y,\cdot)$ is 1, for any $y\in\lin 0, N\rin$.
It comes from the fact that $a_0=1$ and that \bq
\fo n\in\lin 1, N\rin,\qquad \int h_n\, d\gamma&=&0\eq
\wwtbp
\par
Here are two observations about $\cA_N$:
\begin{lem}\label{procAN}
Consider $a\df(a_n)_{n\in\lin 0, N\rin}\in\cA_N$. Then for any
$ n\in\lin 0, N\rin$, $a_n=0$ if $n$ is odd and $a_n\geq 0$ if $n$ is even.
Furthermore, for any $M\in\lin 0, N\rin$, we have $(a_n)_{n\in\lin 0, M\rin}\in\cA_M$.
\end{lem}
\proof
For $y\in\lin 0, N\rin$ odd, $\lambda_a(y,x)$ cannot stay non-negative as $x$ goes to $\pm\iy$ if the coefficient of higher order $a_y$ is non-zero, so we must have $a_y=0$.
For $y\in\lin 0, N\rin$ even, $\lambda_a(y,x)$ cannot stay non-negative if the coefficient of higher order $a_y$ is negative, so we must have $a_y\geq 0$.\par\sm
The second assertion is an immediate consequence of the fact that for $a\in\RR^{N+1}$, $y\in\lin 0, N\rin$ and $x\in\RR$, the definition of the quantity $\lambda_a(y,x)$ does not depend on $N$.
\wwtbp
\par
It follows from the second part of Lemma \ref{procAN} that for any $M\leq N\in\ZZ_+$,
 an element of $\cL_{N}$ can be seen as a element of $\cL_{M}$  by restriction, namely
we have
\bqn{compat}
\Lambda\in\cL_{N}&\Longrightarrow& I_{M, N}\Lambda\in \cL_{M}\eqn
where we recall that $I_{M,N}$ is the natural imbedding from  $\lin 0, M\rin$ into $\lin 0, N\rin$.
To see that $\cL_{M}$ can be strictly larger than $I_{M,N} \cL_{N}$,
let us compute the first sets $\cL_N$. It will appear that $\cL_2$ is strictly included into $I_{2,3} \cL_{3}$.
\par
From Proposition \ref{comput} and Lemma \ref{procAN}, we get that
 $\cL_0= \{\gamma\}$ (the only possible motion on a singleton is to stay still, so an intertwining relation from a singleton
 state space is equivalent to the existence of an invariant measure for the second Markov generator)
 and $\cL_1=\{\gamma\}$.
The next two cases $\cL_2$ and $\cL_3$ are no longer singletons:
\begin{lem}
We have
\bq
\cL_2&=&\Big\{\Lambda\in\cK_2\stl  \Lambda(0,\cdot)= \Lambda(1,\cdot)=\gamma,\, \Lambda(2,dx) =(1+a_2h_2(x)/2)\gamma(dx)\hbox{, with $a_2\in [0, 2]$}\Big\}\eq
and
\bq
\cL_3&=&\Big\{\Lambda\in\cK_3\stl  \Lambda(0,\cdot)= \Lambda(1,\cdot)=\gamma,\, \Lambda(2,dx) =(1+a_2h_2(x)/2)\gamma(dx),\\
&&\Lambda(3,dx) =(1+3a_2h_2(x)/2)\gamma(dx)
\hbox{, with $a_2\in [0,2/3]$}\Big\}\eq
\end{lem}
\proof
From Proposition \ref{comput} and Lemma \ref{procAN}, $\cL_2$ and $\cL_3$ must have the above form, except that for $\cL_2$ (respectively $\cL_3$) the condition on $a_2\geq 0$
is that $\RR\ni x\mapsto 1+a_2h_2(x)/2$ (resp.\ $\RR\ni x\mapsto 1+a_2h_2(x)/2$ and $\RR\ni x\mapsto 1+3a_2h_2(x)/2$) must stay non-negative.\par
Since $h_2(x)=x^2-1$ for all $x\in\RR$, these conditions end up being equivalent to $a_2\in[0,2]$ for $\cL_2$ and
$a_2\in[ 0, 2/3]$ for $\cL_3$.\wwtbp
\par
More generally, consider $N\in \ZZ_+$ with $N\geq 2$.
Since each of the mappings $h_{2n}$, for $n\in\ZZ_+$, is bounded below, we see
that we can find
$(a_{2n})_{n\in\lin 1, \lfloor N/2\rfloor\rin}$ small enough so that the corresponding Markov kernel belongs to $\cL_N$. In particular, $\cL_N$ is not reduced to a singleton, as announced
in (e) of the introduction.\par
\sm
There is no difficulty in replacing the finite sets $\lin 0, N\rin$, for $N\in\ZZ_+$, by $\ZZ_+$. The only supplementary information we need is the following result, extending
Lemma \ref{varphin}:
\begin{lem}\label{specdecomp}
There exist $l\in\CC$ and $f\in\CC^{\ZZ_+}\setminus\{0\}$ such that $D_\iy[f]=-lf$ if and only if $l\in\ZZ_+$
and $f$ is proportional to $\varphi_{l}\in\CC^{\ZZ_+}$ that we recall was  defined  in \eqref{defphi} as 
\bq
\fo y\in\ZZ_+,\qquad \varphi_{l}(y)&\df& \binom{y}{l}\eq
\end{lem}
\proof
Since $f$ is not vanishing everywhere, there exists $N\in\ZZ_+$ such that $f(N)\neq 0$.
It follows that $I_N[f]\in \CC^{\lin 0, N\rin}\setminus\{0\}$, where $I_N$ is the natural embedding of $\lin 0, N\rin$ into $\ZZ_+$. Since $D_NI_N=I_ND$, we deduce that
\bq
D_N[I_N[f]]&=&-lI_N[f]\eq
namely $-l$ is an eigenvalue of $D_N$, thus $l\in \lin 0, N\rin$.
Furthermore, according to Lemma \ref{varphin}, $f$ is proportional to $\varphi_{l}$ on $\lin 0, N\rin$,  say
$f=z\varphi_{l}$ on $\lin 0, N\rin$ for some $z\in \CC$.
Using the relation $D_\iy[f]=-lf$, we show by an iteration on $n\in\ZZ_+$ that
$f=z\varphi_{l}$ on $\lin 0, N+n\rin$.
It follows that $f=z\varphi_{l}$ on $\ZZ_+$. Conversely, there is no difficulty in checking that $D_\iy[\varphi_l]=-l\varphi_l$ for $l\in\ZZ_+$.
\wwtbp
\par
Next, similarly to the beginning of this subsection,
 we define  $\cK_{\iy}$ and $\cL_{\iy}$ as convex sets of Markov kernels from $\ZZ_+$ to $\RR$.
With the help of Lemma \ref{specdecomp}, the proof of Proposition \ref{comput} is still valid and shows that
\bq
\cL_\iy&=&\{\Lambda_a\stl  a\in\cA_\iy\}\eq
where $\Lambda_a$ and $\cA_\iy$ are defined as above, but with $\lin 0, N\rin$ replaced by $\ZZ_+$ in \eqref{lama}, \eqref{Lama} and \eqref{cAN}.
The problem is that it is not clear that $\cA_\iy$ is not reduced to the singleton $\{(\delta_{0,n})_{n\in\ZZ_+}\}$, which is equivalent to say there is  non-trivial intertwinings from $D_\iy$ to $L$.\par
The  compatibility property \eqref{compat}  leads to the following characterization of $\cL_\iy$:
\bq
\fo \Lambda\in \cK_{\iy},\qquad \Lambda\in\cL_{\iy}&\Longleftrightarrow & \fo N\in\ZZ_+,\, I_N\Lambda\in\cL_N\eq
\par
Note that  this latter property with $N=1$ and thus $\cL_1=\{\gamma\}$ contradicts   the first part of Theorem~3.4 in Biane \cite{MR1459446}. Indeed, thereout the author
provides a Markov kernel $\Lambda\in\cL_\iy$  such that
\bq
 \Lambda(1,dx)&=&h_1^2(x)\, \gamma(dx)\
=\ x^2\, \gamma(dx)
\eq
It can also be checked directly with the Markov kernel $\Lambda$ given in \cite{MR1459446}, that $\Lambda[h_2]$ is not proportional to the function $\varphi_2$
defined in Lemma \ref{specdecomp}.
We are thus left wondering if Conjecture \ref{con1} is true, i.e.\ whether $\cL_\iy$ is reduced to a singleton.
\par
\begin{rem}
Let $(a_{2n})_{n\in\ZZ_+}$ be a sequence of non-negative numbers, with $a_0=1$.
Assume there is only a finite number of elements of this sequence that are non-zero.
Then for the mapping
\bq
\RR\ni x&\mapsto &\sum_{n\in\ZZ_+} \frac{a_{2n}}{(2n)!}\binom{y}{n}h_{2n}(x)\eq
to remain non-negative for all $y\in\ZZ_+$, we must have $a_{2n}=0$ for all $n\in\NN$.
Indeed, suppose on the contrary that there exists $N\in\NN$ with $a_{2N}>0$ and consider $N$ the largest such positive integer.
There exists some $x_0\in\RR$ with $h_{2N}(x_0)<0$ (since $\int h_{2N}\, d\gamma=0$).
The quantity $\sum_{n\in\ZZ_+} \frac{a_{2n}}{(2n)!}\binom{y}{n}h_{2n}(x_0)=\sum_{n\in\lin 0, N\rin} \frac{a_{2n}}{(2n)!}\binom{y}{n}h_{2n}(x_0)$
behaves like $ \frac{a_{2N}}{(2N)!}\binom{y}{N}h_{2n}(x_0)<0$ for large $y\in\ZZ_+$, because $\binom{y}{N}\sim y^N/(N!)$.
This is a contradiction justifying the above assertion.
\par
Thus any sequence $a\in \cA_\iy\setminus\{0\}$
 must have an infinite number of non-zero elements.
\end{rem}
\par

Let us discuss the probabilistic implications of the previous considerations.
Fix $N\in\ZZ_+$. For $\Lambda\in\cL_N$, consider $\fP_N(\Lambda)$ the convex hull generated by $\{\Lambda(y,\cdot)\stl   y\in\lin 0,N\rin\}$.
Choose some $\mu_0\in \fP_N(\Lambda)$ and $m_0$ a probability measure on $\lin 0, N\rin$ such that $m_0\Lambda =\mu_0$.
Let $X\df(X_t)_{t\geq 0}$ be an Ornstein-Uhlenbeck diffusion with initial law $\mu_0$ and let $Y\df(Y_t)_{t\geq 0}$ be a Yule jump process with generator $D_N$ and starting from $m_0$.
From the general theory of Markov intertwinings developed by Diaconis and Fill \cite{MR1071805} (see \cite{MR3634282} for an example of  technical extension to a one-dimensional diffusion context),
it is possible to construct a coupling of $X$ and $Y$ such that \eqref{intertwining} is satisfied.
 As a consequence, the stopping time
 \bq
 \tau&\df&\inf\{t\geq 0\stl  Y_t=0\}\eq
 is a \textbf{strong stationary time} for $X$, namely it is a finite randomized stopping time for $X$ such that $\tau$ and $X_\tau$ are independent and $X_\tau$ is distributed according to $\gamma$.
 \par
 Note that the law of $\tau$ is a mixture of the distributions $\cE(1)\convolution\cE(2)\convolution\cdots\convolution\cE(n)$, for $n\in \lin 0, N\rin$, where $\cE(n)$ stands for the exponential law of parameter $n$
 and $\convolution$ for the convolution.
 In particular the law of $\tau$ is stochastically dominated by $\cE(1)\convolution\cE(2)\convolution\cdots\convolution\cE(N)$.
 Recall that the \textbf{separation} $\fs(\mu,\mu')$ between two probability distributions on a same measurable space is defined by
 \bq
 \fs(\mu,\mu')&\df&\esssup_{\mu'}1-\frac{d\mu}{d\mu'}\eq
 where $d\mu/d\mu'$ is the Radon-Nykodim density of $\mu$ with respect to $\mu'$.
 The above observations (see Diaconis and Fill \cite{MR1071805} for the general argument) lead to the following quantitative estimates on the convergence to equilibrium for the Ornstein-Uhlenbeck diffusion
 in the total variation and separation senses:
 \bq
 \fo t\geq 0,\qquad \lVe \cL(X_t)-\gamma\rVe_{\mathrm{tv}}&\leq & \fs(\cL(X_t),\gamma)\\
 &\leq &\PP[\tau>t]\\
 &\leq & \cE(1)\convolution\cE(2)\convolution\cdots\convolution\cE(N)((t,+\iy))\eq
 under the condition that the initial law $\cL(X_0)$ belongs to
 \bq
 \fP_N&\df& \bigcup_{\Lambda\in\cL_N} \fP_N(\Lambda)\eq
 \par\sm
 Considering
 \bq
 \fP&\df & \bigcup_{n\in\ZZ_+}\fP_N\eq
 we get that when the initial distribution of $X$ belongs to $\fP$, there exists a strong stationary time for $X$.
 It was proven in \cite{MR3634282} this is not true for all initial distributions, in particular for those with compact support.
 Nevertheless, one can  wonder if $\fP$ would not be dense, e.g.\  in the total variation sense, in the set of probability measures on $\RR$ absolutely continuous with respect to $\gamma$.
 Other natural questions are:
 \begin{description}
 \item[$\bullet$] is the set $\fP_N$ convex and in this case what are the extremal points?
 \item[$\bullet$] is the sequence $(\fP_N)_{N\in\ZZ_+}$ non-decreasing?
 \item[$\bullet$] what is the link between $\fP$ and $\fP_\iy\df \bigcup_{\Lambda\in\cL_\iy} \fP_\iy(\Lambda)$?, where for $\Lambda\in\cL_\iy$, $\fP_\iy(\Lambda)$ is the set of $m\Lambda$, with $m$ a probability measure on $\ZZ_+$.
 \end{description}

\subsection{On the reverse Yule family}\label{otrYf}

Under the questions at the end of last subsection lies the interrogation: is it always possible to slightly modify the initial condition
of an Ornstein-Uhlenbeck diffusion to insure the existence of a strong stationary time?
A reverse question is: are there strong times which bring the diffusion close to its equilibrium?
Recall that a \textbf{strong time} is a finite randomized stopping time $\tau$ for $X$ such that
$\tau$ and $X_\tau$ are independent. The reverse Yule family seems more appropriate to construct strong times that deal first
with the low lying eigenvalues and next with the high lying eigenvalues: the corresponding jump processes begin by  encountering the rate 1, then the rate 2, etc.
It was the opposite with the Yule processes of the previous sections.
Thus in the spirit of the motivations described in the introduction, reverse Yule processes should be preferable.
Note this point of view can also be found in Fill \cite{MR2530102} and in \cite{MR2654550}, but in finite settings.
Unfortunately, due to (f) of the introduction, this direction ends up being not so relevant, except if Conjecture \ref{conrev}  was to be true.
\par\me
Let us try to extend the analysis done for the Yule family to the reverse Yule family.
Our first task is to compute the eigenvectors of the reverse Yule family.
\begin{lem}
Fix $N\in\ZZ_+$.
For all $n\in \lin N\rin\df \lin 1, N\rin$, the eigenspace associated to the eigenvalue $-n$ of $\wi D_N$, defined in \eqref{wiDN}, is generated by the function $\wi\varphi_n$ given by
\bq
\fo y\in \lin 0, N\rin,\qquad \wi\varphi_n(y)&\df&
 (-1)^y
\binom{n-1}{-y}\eq
For $n=0$, it is sufficient to take $\wi\varphi_0=\un$.
\end{lem}
\proof
For $n=0$, it is obvious that $\wi D_N[\un]=0$.
For $n\in\lin N\rin$, the relation $\wi D_N[\wi\varphi_n]=-n\wi\varphi_n$ leads to an iteration on the values of $\wi \varphi_n$:
\bq
\fo y\in\lin -(N-1),0\rin,\qquad \wi \varphi_n(y-1)&=&\frac{1-n-y}{1-y}\wi \varphi_n(y)\eq
which gives the announced function, when we take $\wi \varphi_n(0)=1$.
\wwtbp
\par
Remark that $\wi \varphi_n(y)=0$ as soon as $y\leq -n$, in particular all $\wi \varphi_n$, for $n\in\lin N\rin$, vanish at $-N$. This comes from the fact that $-N$ is absorbing for $\wi D_N$.
Due to this property, the restriction of $\wi D_\iy$ to $\lin -N, 0\rin$ is different from $\wi D_N$, but only at the entry $(-N,-N)$.
So the proof by restriction of Lemma~\ref{specdecomp} cannot be applied directly, nevertheless we get a similar  result (either by a direct proof or by taking into account that
$\wi \varphi_n(-N)=0$ for all $n\in\lin N\rin$):
the spectrum of $-\wi D_\iy$ is $\ZZ_+$, $\un$ spans the eigenspace associated to 0 and for $n\in\NN$, the eigenspace associated to $n$ is generated by the extension to $\ZZ_-$
of the previous function $\wi\varphi_n$, namely:
\bq
\fo y\in \ZZ_-,\qquad \wi\varphi_n(y)&\df&
 (-1)^y
\binom{n-1}{-y}\eq\par\me
As in the beginning of Subsection \ref{fYtO}, consider $\wi\cK_N$
the set of Markov kernels from $\lin -N,0\rin$ to $\RR$ such that for all $y\in \lin -N,0\rin$, $\Lambda(y,\cdot)\in \cM$, and
\bq
\wi\cL_N&\df&\{ \Lambda\in\wi\cK_N\stl  \wi D_N\Lambda=\Lambda L\}\eq
Introduce
for any $a\df(a_n)_{n\in\lin 0,N\rin}$,  the mapping $\wi\lambda_a\stl  \lin0,N\rin\times\RR\ri\RR$
and the signed kernel $\wi \Lambda_a$ from $\lin -N, 0\rin$ to $\RR$ via
\bqn{wilama}
\fo y\in\lin 0, N\rin,\,\fo x\in \RR,\qquad \wi\lambda_a(y,x)&\df& \sum_{n\in\lin -y+1,N\rin} \frac{a_n}{n!}\wi \varphi_n(y)h_n(x)\\
\label{wiLama} \fo y\in\lin 0, N\rin,\qquad \wi\Lambda_a(y,dx)&\df& \wi\lambda_a(y,x)\, \gamma(dx)
\eqn
Finally, write $\wi\cA_N$ for the set of $a\df(a_n)_{n\in\lin 0,N\rin}\in \RR^{N+1}$ with $a_0=1$ such that
\bqn{wicAN}
\fo y\in\lin 0, N\rin,\,\fo x\in \RR,\qquad \wi\lambda_a(y,x)&\geq &0\eqn
 \par
The proof of Proposition \ref{comput} leads to
\bq
\wi\cL_N&=&\{\wi\Lambda_a \stl  a\in\wi\cA_N\}\eq
\par
Thus Assertion (f) of the introduction is a consequence of:
\begin{lem}
We have
\bq
\wi\cA_N&=&\{(1, 0, 0, ..., 0)\}
\eq
\end{lem}
\proof
Fix $a\df(a_n)_{n\in\lin 0,N\rin}\in\wi\cA_N$ and consider
\bq
n_0&\df&\max\{n\in\lin 0, N\rin\stl  a_n\neq 0\}\eq
Assume $n_0\geq 1$ and let us obtain a contradiction.
For $y\in\lin -n_0+1, 0\rin$, we have $\wi\varphi_{n_0}(y)\neq 0$,
thus as $x$ goes to $\pm\iy$, $\lambda_a(y,x)$ is equivalent to $a_{n_0}\wi\varphi_{n_0}(y)x^{n_0}/n_0!$
and we must have $n_0$ even and $\wi\varphi_{n_0}(y)>0$.
It follows that $n_0\geq 2$ and then $\lin -n_0+1,0\rin$ contains odd and even points $y$.
Since the sign of $\wi\varphi_{n_0}(y)$ is $(-1)^y$, we end up with a contradiction.\wwtbp
\par
\sm
The case $N=\iy$ is not so clear, since now \eqref{wilama} is an infinite sum and we don't know precisely in which sense it should converge,
probably in some weak sense so that \eqref{wiLama} has still a meaning (maybe with the $\wi\Lambda_a(y,\cdot)$ no longer absolutely continuous with respect $\gamma$).
The set $\wi\cA_\iy$ should be the set of all $a\df(a_n)_{n\in\ZZ_+}\in \RR^{\ZZ_+}$ with $a_0=1$ such that
$\wi\Lambda_a$ is a non-negative kernel and we would end up with the conclusion that
$\wi\cL_\iy=\{\wi\Lambda_a \stl  a\in\wi\cA_\iy\}$.
The above arguments are no longer sufficient to deduce that $\wi\cL_\iy$ is a singleton.
On the contrary,
we believe that Conjecture \ref{conrev} is true.


\vskip2cm
\hskip60mm
\vbox{
\copy5
 \vskip5mm
 \copy6
}

\end{document}